\documentclass[twoside,11pt]{article}
\title{} \author{} \date{}
\usepackage{amssymb,amsmath,latexsym,times}%,showkeys}
\newtheorem{te}{Theorem}[section]
\newtheorem{fac}[te]{Fact}

\newtheorem{ex}[te]{Example}

\def\upar{\!\uparrow}
\def\downar{\!\downarrow}
\def\dok{\noindent{\bf Proof. }}
\def\kdok{\hfill $\Box$ \par \vspace*{2mm} }
\def\Lim{\lim\nolimits}
\newcommand{\Conv}{\mathop{\mathrm{Conv}}\nolimits}
\newcommand{\TopConv}{\mathop{\mathrm{TopConv}}\nolimits}
\newcommand{\Top}{\mathop{\mathrm{Top}}\nolimits}
\newcommand{\SeqTop}{\mathop{\mathrm{SeqTop}}\nolimits}
\newcommand{\Br}{\mathop{{\mathrm{Br}}}}
\newcommand{\Coll}{\mathop{{\mathrm{Coll}}}}
\begin{document}
\thispagestyle{plain}
\begin{center}
           {\large \bf
           \uppercase{The left, the right and the sequential topology \\[1mm] on Boolean  algebras}}
\end{center}
\begin{center}
{\bf Milo\v s S.\ Kurili\'c\footnote{Department of Mathematics and Informatics, Faculty of Science, University of Novi Sad,
Trg Dositeja Obradovi\'ca 4, 21000 Novi Sad, Serbia. e-mail: milos@dmi.uns.ac.rs}
and Aleksandar Pavlovi\'c\footnote{Department of Mathematics and Informatics, Faculty of Science, University of Novi Sad,
Trg Dositeja Obradovi\'ca 4, 21000 Novi Sad, Serbia. e-mail: apavlovic@dmi.uns.ac.rs}}
\end{center}
\begin{abstract}
\noindent
For the algebraic convergence $\lambda_{\mathrm{s}}$,
which generates the well known sequential topology $\tau_s$ on a complete Boolean algebra ${\mathbb B}$,
we have $\lambda_{\mathrm{s}}=\lambda_{\mathrm{ls}}\cap \lambda_{\mathrm{li}}$, where the convergences
$\lambda_{\mathrm{ls}}$ and $\lambda_{\mathrm{li}}$ are defined by
$\lambda_{\mathrm{ls}}(x)=\{ \limsup x\}\!\upar$ and $\lambda_{\mathrm{li}}(x)=\{ \liminf x\}\downar$
(generalizing the convergence of sequences on the Alexandrov cube and its dual).
We consider the minimal topology $\mathcal{O}_{\mathrm{lsi}}$ extending
the (unique) sequential topologies $\mathcal{O}_{\lambda_{\mathrm{ls}}}$ (left) and $\mathcal{O}_{\lambda_{\mathrm{li}}}$ (right)
generated by the convergences $\lambda_{\mathrm{ls}}$ and $\lambda_{\mathrm{li}}$
and establish a general hierarchy between all these topologies and the corresponding a priori and a posteriori convergences.
In addition, we observe some special classes of algebras and, in particular, show that
in $(\omega,2)$-distributive algebras we have $\lim_{{\mathcal O}_{\mathrm{lsi}}}=\lim_{\tau _{\mathrm{s}} }=\lambda _{\mathrm{s}}$,
while the equality $\mathcal{O}_{\mathrm{lsi}}=\tau_s$ holds in all Maharam algebras. On the other hand,
in some collapsing algebras we have a maximal (possible) diversity. 
%\\[1mm]

%\\%[1mm]

\end{abstract}
\section{Introduction}
It is known
%(and easy to check)
that a sequence $\langle x_n :n\in \omega\rangle$ of reals from the unit interval $I=[0,1]$
converges to a point $a\in I$ with respect to the left (resp.\ right, standard) topology on $I$ iff
$a\geq \limsup x_n$ (resp.\ $a\leq \liminf x_n$, $a= \liminf x_n =\limsup x_n$) and, more generally, these three properties
define three convergence structures on any complete lattice or $\sigma$-complete Boolean algebra.
In this paper, continuing the investigation from \cite{KuPa07}--\cite{KuTo}, we consider the corresponding
convergences $\lambda_{\mathrm{ls}}$, $\lambda_{\mathrm{li}}$ and $\lambda_{\mathrm{s}}$
on a  complete Boolean algebra $\mathbb{B}$, as well as the sequential topologies
$\mathcal{O}_{\lambda_{\mathrm{ls}}}$, $\mathcal{O}_{\lambda_{\mathrm{li}}}$ and $\mathcal{O}_{\lambda_{\mathrm{s}}}$ on $\mathbb{B}$
generated by them. Having in mind that the union of the left and the right topology on $I$ generates the standard topology on that interval,
we regard the minimal topology $\mathcal{O}_{{\mathrm{lsi}}}$ on $\mathbb{B}$ extending
$\mathcal{O}_{\lambda_{\mathrm{ls}}} \cup \mathcal{O}_{\lambda_{\mathrm{li}}}$, as well as the corresponding
topological convergence $\lim_{\mathcal{O}_{{\mathrm{lsi}}}}$ on $\mathbb{B}$, and explore the relationship between all the
topologies and convergences mentioned above. It turns out that everything consistent is possible.
For example, if $\mathbb{B}$ is the power set algebra $P(\omega )$ or a Maharam algebra (i.e.\ admits a strictly positive Maharam submeasure),
then we have an analogy to the unit interval: $\mathcal{O}_{{\mathrm{lsi}}}=\mathcal{O}_{\lambda_{\mathrm{s}}}$; but for some collapsing algebras we obtain a
maximal diversity (e.g.\ $\lim_{\mathcal{O}_{{\mathrm{lsi}}}}\neq \lim_{\mathcal{O}_{\lambda_{\mathrm{s}}}}$ and % , hence,
$\mathcal{O}_{{\mathrm{lsi}}}\neq \mathcal{O}_{\lambda_{\mathrm{s}}}$).

We note that the  topology $\mathcal{O}_{\lambda_{\mathrm{s}}}$ (traditionally called the {\it sequential topology on c.B.a,'s} and denoted by $\tau_s$),
generated by the convergence $\lambda_{\mathrm{s}}$ (traditionally called the {\it algebraic convergence})
was widely considered in the context of the von Neumann problem \cite{Sco81}: Is each ccc  weakly distributive c.B.a.\  a measure algebra?
A consistent counter-example (a Suslin algebra) was given by  Maharam  \cite{Mah47}.
In addition, Maharam has shown that the topology $\mathcal{O}_{\lambda_{\mathrm{s}}}$ is metrizable
iff $\mathbb{B}$ is a Maharam algebra and asked
whether this implies that  $\mathbb{B}$ admits a measure (the Control Measure Problem, negatively solved by M.\ Talagrand \cite{Tal1, Tal2}).
Moreover, Balcar, Jech and Paz\'{a}k \cite{BJP05} and, independently, Veli\v ckovi\'c \cite{Vel05}, proved that it is consistent that
the topology $\mathcal{O}_{\lambda_{\mathrm{s}}}$ is metrizable on each complete ccc weakly distributive algebra.
(See also \cite{BGJ98,BJ06,Far04,Tod04} for that topic).

Regarding the power set algebras, $P(\kappa)$, the convergence $\lambda_{\mathrm{s}}$ is exactly the convergence on the Cantor cube,
while $\lambda_{\mathrm{ls}}$ generalizes the convergence on the Alexandrov cube in the same way (see \cite{KuPaAlx}).
Further, on any c.B.a., the topologies ${\mathcal O}_{\lambda_{\mathrm{ls} }}$ and ${\mathcal O}_{\lambda_{\mathrm{li} }}$
are homeomorphic (take $f(a)=a'$) and  generated by
some other convergences relevant for set-theoretic forcing (see \cite{KuPaNSJOM,KuPaDeb}).
For obvious reasons,  the topology ${\mathcal O}_{\lambda_{\mathrm{ls} }}$ (resp.\ ${\mathcal O}_{\lambda_{\mathrm{li} }}$)
will be called the {\it left} (resp.\ the {\it right}) {\it topology on} $\mathbb{B}$ (see also Fact \ref{T1261}(i)).
\section{Preliminaries}
\paragraph{Convergence}
Here we list the standard facts concerning convergence structures which will be used in the paper.
(For details and proofs see, for example, \cite{KuPaNSJOM}.)

Let $X$ be a non-empty set. Each mapping $x:\omega \to X$ is called a {\it sequence} in $X$.
Usually, instead of $x(n)$ we write $x_n$ and $x=\langle x_n : n\in \omega\rangle$.
A {\it constant sequence} $\langle a,a,\ldots \rangle$ is denoted shortly by $\langle a\rangle$.
A sequence $y \in X^\omega$ is said to be a {\it subsequence} of $x$ iff there is an increasing function $f: \omega\rightarrow \omega$
(notation: $f\in \omega^{\uparrow \omega}$)
such that $y=x \circ f$; then we write $y \prec x$.

Each mapping $\lambda: X^\omega \to P(X)$ is called a {\it convergence}. The set $\Conv (X)=P(X)^{(X^\omega)}$ of all convergences on the set $X$
ordered by the relation $\lambda _1\leq \lambda_2$ iff $\lambda _1(x)\subseteq \lambda_2(x)$, for each $x\in X^\omega$, is, clearly, a Boolean lattice
and $\lambda _1 \cap \lambda _2 $ will denote
the infimum $\lambda _1 \wedge \lambda _2 $; that is, $(\lambda _1 \cap \lambda _2) (x) =\lambda _1 (x) \cap \lambda _2 (x)$, for all $x\in X^\omega$.
If $|\lambda(x) |\leq 1$ for each sequence $x$, then $\lambda$ is called a {\it   Hausdorff convergence}.

Let $\langle X, {\mathcal O}\rangle$ be a topological space.
A point $a \in X$ is a {\it limit point} of a sequence $x\in X^\omega$ iff
each neighborhood of $a$ contains all but finitely many members of $x$.
The set of all limit points of a sequence $x\in X^\omega$ is denoted by $\lim _{\mathcal O}(x)$
and so we obtain a convergence $\lim _{\mathcal O}: X^\omega \rightarrow P(X)$, that is, $\lim _{\mathcal O}\in \Conv (X)$.

Let $\Top (X)$ denote the lattice of all topologies on the set $X$.
A convergence $\lambda\in \Conv (X)$ is called {\it topological}, we will write $\lambda \in \TopConv (X)$, iff
there is a topology ${\mathcal O}\in \Top (X)$ such that $\lambda=\lim_{\mathcal O}$. So we establish the mapping
$$\textstyle
G: \Top (X) \rightarrow  \TopConv (X), \mbox{ where } G(\mathcal O )=\lim _{\mathcal O}.
$$
A topology ${\mathcal O}\in\Top (X)$ is called {\it sequential}, we will write ${\mathcal O}\in\SeqTop (X)$
iff in the space $\langle X, {\mathcal O}\rangle$ we have:
a set $A\subset X$ is closed iff  it is {\it sequentially closed} (that is, $\lim _{\mathcal O}(x) \subset A$, for each sequence $x\in A^\omega$).
If ${\mathcal O}_1, {\mathcal O}_2 \in\SeqTop (X)$ and $\lim_{{\mathcal O}_1} = \lim_{{\mathcal O}_2}$, then ${\mathcal O}_1={\mathcal O}_2$.\footnote{
This is false in general: take  the discrete and the co-countable topology on the real line; in both spaces exactly the almost-constant sequences converge.}
So, $G$ is one-to-one on $\SeqTop (X)$.

For each convergence $\lambda\in \Conv (X)$ there is a (unique) maximal topology ${\mathcal O}_\lambda$ such that
$\lambda \leq \lim_{{\mathcal O}_\lambda}$. The topology $\mathcal{O}_\lambda $ is sequential; so, we obtain the mapping
$$
F: \Conv (X)\rightarrow \SeqTop (X), \mbox{ defined by } F(\lambda )=\mathcal{O}_\lambda .
$$
$F$ and $G$ are antitone mappings, that is, $\lambda _1 \leq \lambda _2$ implies that
%\Rightarrow
$\mathcal{O}_{\lambda _2} \subset \mathcal{O}_{\lambda _1}$
and $\mathcal{O}_1 \subset \mathcal{O}_2$
%\Rightarrow
implies $\lim _{{\mathcal O}_2} \leq \lim _{{\mathcal O}_1} $.
Moreover, a convergence $\lambda $ is topological iff $\lambda =\lim _{{\mathcal O}_\lambda}(=G(F(\lambda )))$ and
a topology ${\mathcal O}$ is sequential iff ${\mathcal O}={\mathcal O}_{\lim _{\mathcal O}}(= F(G({\mathcal O}))$\footnote{In fact,
the pair $F,G$ is an antitone Galois connection between the complete lattices $\Conv (X)$ and $\Top (X)$,
because
$\mathcal{O} \subset F(\lambda ) \Leftrightarrow \lambda \leq G({\mathcal O})$, for each $\lambda \in \Conv (X)$ and ${\mathcal O}\in \Top (X)$.
(If $\mathcal{O} \subset {\mathcal O}_\lambda $, then
$\lambda \leq \lim_{{\mathcal O}_\lambda}\leq \lim_{{\mathcal O}}$.
Conversely, if $\lambda \leq \lim_{{\mathcal O}}$, then $\mathcal{O} \subset {\mathcal O}_\lambda $, by the maximality of ${\mathcal O}_\lambda $).
Moreover, the restriction $F\upharpoonright \TopConv (X)$ is a bijection from $\TopConv (X)$ onto $\SeqTop (X)$
and $G\upharpoonright \SeqTop (X)$ is its inverse.
}.
Each topological convergence $\lambda$ satisfies the following conditions:

(L1) $\forall a \in X\;\;  a \in \lambda (\langle a\rangle )$,

(L2) $\forall x \in X^\omega\; \forall y \prec x\; \lambda (x) \subset \lambda (y)$,

(L3) $\forall x \in X^\omega \;
     \forall a \in X\;
     ((\forall y \prec x\;
     \exists z\prec y \;
     a \in \lambda(z)) \Rightarrow a \in \lambda(x))$.

\noindent
If $\lambda \in \Conv (X)$ satisfies (L1) and (L2),
then $\mathcal{O}_\lambda \!=\!\{X\setminus F: F\subset X \land u_\lambda(F)\!=\!F\}$, where $u_\lambda:P(X) \to P(X)$
is the {\it operator of sequential closure determined by $\lambda$}, defined by $u_\lambda(A)=\bigcup_{x\in A^\omega}\lambda(x)$.
In addition, the minimal closure of $\lambda$ under (L1)--(L3) is given by
$\textstyle \lambda^*(x)=\bigcap_{f\in \omega^{\uparrow \omega}}\bigcup_{g\in \omega^{\uparrow \omega}}{\lambda}(x \circ f \circ g)$ and
$\lambda$ is called a {\it  weakly-topological convergence} iff the convergence $ \lambda^*$ is topological.
\begin{fac}\label{T1220}
(\cite{KuPaNSJOM})
If $\lambda \in \Conv (X)$ is a convergence satisfying (L1) and (L2), then

(a) $\lambda$ is weakly-topological iff $\,\lim_{{\mathcal O}_\lambda }= \lambda^*$, that is, for each $x\in X^\omega$ and $a\in X$
$$
\textstyle a \in \lim_{{\mathcal O}_\lambda}(x) \Leftrightarrow
\forall y \prec x \;\; \exists z \prec y \;\; a \in \lambda(z);
$$

(b) If $\lambda$ is a Hausdorff convergence, then $\lambda^*$ is Hausdorff and weakly-topological.
\end{fac}
\paragraph{Convergences on Boolean algebras} Let ${\mathbb B}$ be a complete Boolean algebra or, more generally,
a complete lattice. If $\langle x_n :n\in \omega \rangle$ is a sequence of its elements,
$\liminf x_n  := \textstyle \bigvee _{k\in \omega}\bigwedge _{n\geq k}x_n$ and
$\limsup x_n  := \textstyle\bigwedge _{k\in \omega}\bigvee _{n\geq k}x_n $, then, clearly,  $\liminf x_n \leq \limsup x_n$.
We consider the convergences $\lambda_{\mathrm{\mathrm{ls}} }, \lambda_{\mathrm{\mathrm{li}} }, \lambda_{\mathrm{s} }:B^\omega \rightarrow P(B)$ defined by
\begin{eqnarray}
\lambda_{\mathrm{ls} }(\langle x_n \rangle) & = & \{ \limsup x_n \}\upar ,\\
\lambda_{\mathrm{li} }(\langle x_n \rangle) & = & \{ \liminf x_n \} \downar , \\
\lambda_{\mathrm s}(\langle x_n \rangle)& = & \left\{
                \begin{array}{cl}
                  \{ x \} & \mbox{ if }{\liminf x_n =\limsup x_n =x ,} \\
                  0       & \mbox{ if } {\liminf x_n <\limsup x_n ,}
                \end{array}
              \right.
\end{eqnarray}
where  $A\upar :=\{b\in \mathbb{B}: \exists a \in A \; b \geq a\}$ and
$A\downar :=\{b\in \mathbb{B}: \exists a \in A \; b \leq a\}$, for $A\subset {\mathbb B}$.
The following property of c.B.a.'s will play a role in this paper\footnote{
We note that property ($\hbar$) is closely related to the cellularity of Boolean
algebras. Namely, by \cite{KuPa07}, $\mathfrak{t}\mbox{-cc} \Rightarrow (\hbar ) \Rightarrow \mathfrak{s}\mbox{-cc}$
and, in particular, ccc complete Boolean algebras satisfy ($\hbar$).
By \cite{KuTo}, the set $\{ \kappa \in {\mathrm{Card}} :  \kappa \mbox{-cc }\Rightarrow (\hbar ) \}$
is equal either to $[0, {\mathfrak h})$,  or to $ [0,{\mathfrak h}]$
and $\{ \kappa \in {\mathrm{Card}} : (\hbar )\Rightarrow \kappa \mbox{-cc }    \} =[{\mathfrak s}, \infty )$.
Basic facts concerning the invariants of the continuum
$\mathfrak{t}$, $\mathfrak{s}$, and $\mathfrak{h}$ can be found in \cite{Douw84}.
}
$$
\forall x \in {\mathbb B}^\omega ~ \exists y \prec x~ \forall z \prec y ~ \limsup z =\limsup y .
\eqno{( \hbar )}
$$
\begin{fac} \label{T1212}
(\cite{KuPa07}) If $\mathbb{B}$ is a complete Boolean algebra, then  we have

(a) $\lambda_{\mathrm{s}}$ is a weakly-topological Hausdorff convergence satisfying (L1) and (L2);

(b) $\lambda_{\mathrm{s}}$ is a topological convergence iff the algebra $\mathbb{B}$ is $(\omega,2)$-distributive.
\end{fac}
\begin{fac} \label{T1261}(\cite{KuPaAlx})
If ${\mathbb B}$ is a complete non-trivial Boolean algebra, then

(a) $\lambda _{\mathrm{ls} }$ and $\lambda _{\mathrm{li} }$ are non-Hausdorff convergences satisfying (L1) and (L2);

(b) If $\,{\mathbb B}$ satisfies ($\hbar$), then $\lambda _{\mathrm{ls} }$ and $\lambda _{\mathrm{li} }$ are weakly-topological convergences;

%(c) $\lambda_{\mathrm{ls}}$ is a topological convergence iff $\lambda_{\mathrm{li}}$ is a topological convergence iff
%    the algebra ${\mathbb B}$ is $(\omega,2)$-distributive;

(c) $\lambda_{\mathrm{ls}}$ is topological iff $\lambda_{\mathrm{li}}$ is topological iff the algebra
    ${\mathbb B}$ is $(\omega,2)$-distributive;

(d) $\lambda_{\mathrm s} = \lambda_{\mathrm{ls} } \cap \lambda_{\mathrm{li} }$;

(e) ${\mathcal O}_{\lambda_{\mathrm{ls} }}, {\mathcal O}_{\lambda_{\mathrm{li} }}
     \subset {\mathcal O}_{\lambda_{\mathrm s}}$;

(f) $\lambda _{\mathrm{ls} }^* \leq \lim _{{\mathcal O}_{\lambda _{\mathrm{ls} }}}$ and
    $\lambda _{\mathrm{li} }^* \leq \lim _{{\mathcal O}_{\lambda _{\mathrm{li} }}}$;

(g) $\lambda_{\mathrm s}^* = \lambda_{\mathrm{ls} }^* \cap \lambda_{\mathrm{li} }^*$;

(h)  $\mathcal{O}_{\lambda_{\mathrm{ls}}}$ and $\mathcal{O}_{\lambda_{\mathrm{li}}}$ are homeomorphic, $T_0$, connected and compact  topologies;

(i) A set $F\subset {\mathbb B}$ is ${\mathcal O}_{\lambda_{\mathrm{ls} }}$-closed  iff it is upward-closed and $\bigwedge _{n\in \omega}x_n \in F$, for each decreasing sequence $\langle x_n \rangle\in F^\omega$; (and dually, for ${\mathcal O}_{\lambda_{\mathrm{li} }}$-closed sets).
\end{fac}
\section{The topology ${\mathcal O}_{\mathrm{lsi}}$ on Boolean algebras}
On a complete Boolean algebra ${\mathbb B}$ we consider the minimal topology containing the topologies
${\mathcal O}_{\lambda_{\mathrm{ls} }}$ and ${\mathcal O}_{\lambda_{\mathrm{li}}}$. This topology, denoted by
${\mathcal O}_{\mathrm{lsi}}$, is generated by the  base
${\mathcal B}_{\mathrm{lsi}}= \{ O_1 \cap O_2 : O_1 \in {\mathcal O}_{\lambda_{\mathrm{ls} }}\land O_2 \in {\mathcal O}_{\lambda_{\mathrm{li} }}\}$.
By Fact \ref{T1261}(i), the sets from ${\mathcal O}_{\lambda_{\mathrm{ls} }}$ (resp.\ ${\mathcal O}_{\lambda_{\mathrm{li} }}$)
are downward (resp.\ upward)-closed; so, the elements of  ${\mathcal B}_{\mathrm{lsi}}$ are convex subsets of ${\mathbb B}$.
\begin{te}\label{T1286}
The following diagrams show the relations between the considered convergences and topologies on a non-trivial c.B.a.\ ${\mathbb B}$.
In addition, we have

(a) $\lambda _{\mathrm{ls}} \cap \lambda _{\mathrm{li}} = \lambda _{\mathrm{s}}$,
$\lambda _{\mathrm{ls}} ^*\cap \lambda _{\mathrm{li}}^* = \lambda _{\mathrm{s}}^*$ and
$\lim_{{\mathcal O}_{\lambda _{\mathrm{ls}}}} \cap \lim_{{\mathcal O}_{\lambda _{\mathrm{li}}}}=\lim_{{\mathcal O}_{\mathrm{lsi}} }$;

(b) $\lambda_{\mathrm s}< \lambda_{\mathrm{ls}} , \lambda_{\mathrm{li}}$, $\lambda_{\mathrm s}^*< \lambda_{\mathrm{ls}}^* , \lambda_{\mathrm{li}}^*$,
$\lim_{{\mathcal O}_{\mathrm{lsi}}} < \lim_{{\mathcal O}_{\lambda_{\mathrm{ls} }}} , \lim_{{\mathcal O}_{\lambda_{\mathrm{li} }}}$ and
${\mathcal O}_{\mathrm{lsi}} \supsetneq {\mathcal O}_{\lambda_{\mathrm{ls} }} , {\mathcal O}_{\lambda_{\mathrm{li} }}$.
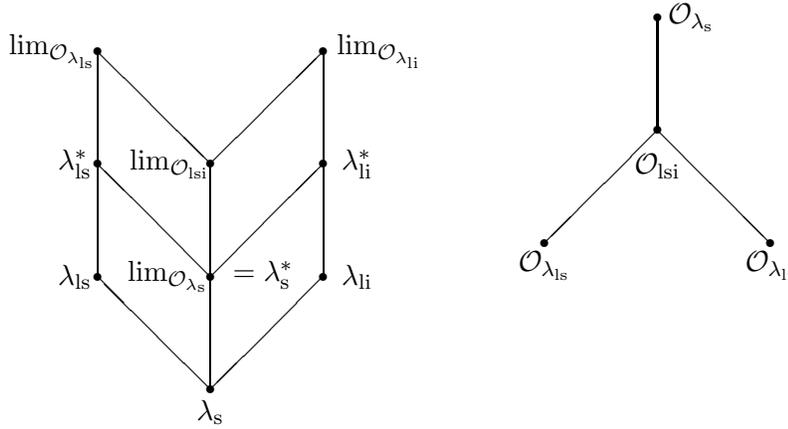
\begin{figure}[h]
  \begin{center}
%==============================================================================
\unitlength .3mm % = .854pt
\linethickness{0.4pt}
\ifx\plotpoint\undefined\newsavebox{\plotpoint}\fi % GNUPLOT compatibility
\begin{picture}(175,276)(0,0)
%===================== linije =================================
\put(150,275){\line(0,-1){50}}
\put(150,225){\line(0,-1){50}}
\put(150,175){\line(-1,-1){50}}
\put(100,125){\line(-1,1){50}}
\put(50,175){\line(0,1){100}}
\put(50,275){\line(1,-1){50}}
\put(100,225){\line(1,1){50}}
\put(100,125){\line(0,1){100}}
\put(100,175){\line(1,1){50}}
\put(100,175){\line(-1,1){50}}
%===================== tacke =================================
\put(100,225){\circle*{4}}
\put(150,275){\circle*{4}}
\put(150,225){\circle*{4}}
\put(150,175){\circle*{4}}
\put(100,125){\circle*{4}}
\put(100,175){\circle*{4}}
\put(50,225){\circle*{4}}
\put(50,275){\circle*{4}}
\put(50,175){\circle*{4}}
%===================== tekst =================================
\put(165,225){\makebox(0,0)[cc]{$\lambda^*_{\mathrm{li}}$}}
\put(40,225){\makebox(0,0)[cc]{${\lambda^*_{\mathrm{ls}}}$}}
\put(100,115){\makebox(0,0)[cc]{$\lambda_{\mathrm{s}}$}}
\put(40,175){\makebox(0,0)[cc]{$\lambda_{\mathrm{ls}}$}}
\put(165,175){\makebox(0,0)[cc]{$\lambda_{\mathrm{li}}$}}
\put(100,175){\makebox(0,0)[cc]{$\lim_{{\mathcal O}_{\lambda_{\mathrm{s}}}}\;\;=\lambda^*_{\mathrm{s}}$}}
\put(82,225){\makebox(0,0)[cc]{$\lim_{{\mathcal O}_{\mathrm{lsi}}}$}}
\put(175,275){\makebox(0,0)[cc]{$\lim_{{\mathcal O}_{\lambda_{\mathrm{li}}}}$}}
\put(30,275){\makebox(0,0)[cc]{$\lim_{{\mathcal O}_{\lambda_{\mathrm{ls}}}}$}}
\end{picture}
      \hspace{5mm}
%================================================================================
%\unitlength .3mm % = .854pt
%\linethickness{0.4pt}
\ifx\plotpoint\undefined\newsavebox{\plotpoint}\fi % GNUPLOT compatibility
\begin{picture}(151,291)(0,0)
%===================== linije =================================
\put(100,290){\line(0,-1){50}}
\put(100,240){\line(-1,-1){50}}
\put(100,240){\line(1,-1){50}}
%===================== tacke =================================
\put(150,190){\circle*{4}}
\put(50,190){\circle*{4}}
\put(100,240){\circle*{4}}
\put(100,290){\circle*{4}}
%===================== tekst =================================
\put(150,180){\makebox(0,0)[cc]{${\mathcal O}_{\lambda_{\mathrm{li}}}$}}
\put(50,180){\makebox(0,0)[cc]{${\mathcal O}_{\lambda_{\mathrm{ls}}}$}}
\put(100,224){\makebox(0,0)[cc]{${\mathcal O}_{\mathrm{lsi}}$}}
\put(115,290){\makebox(0,0)[cc]{${\mathcal O}_{\lambda_{\mathrm{s}}}$}}
\end{picture}
  \end{center}
\vspace{-35mm}
  \caption{Convergences and topologies on ${\mathbb B}$}\label{FIG1}
\end{figure}
\end{te}
\dok
By Fact \ref{T1261}(e) we have
${\mathcal O}_{\lambda_{\mathrm{ls} }}, {\mathcal O}_{\lambda_{\mathrm{li} }}  \subset {\mathcal O}_{\lambda_{\mathrm s}}$
and the inclusion ${\mathcal O}_{\mathrm{lsi}} \subset {\mathcal O}_{\lambda_{\mathrm s}}$ follows from the minimality of ${\mathcal O}_{\mathrm{lsi}} $.
So the diagram for topologies is correct.

By Fact \ref{T1261}(d) and (g) we have
$\lambda _{\mathrm{s}} = \lambda _{\mathrm{ls}} \cap \lambda_{\mathrm{li}} $ and
$\lambda _{\mathrm{s}}^* = \lambda _{\mathrm{ls}}^* \cap \lambda _{\mathrm{li}}^*$,
which implies
$\lambda _{\mathrm{s}}\leq \lambda _{\mathrm{ls}},\lambda_{\mathrm{li}}$
and $\lambda _{\mathrm{s}}^* \leq \lambda _{\mathrm{ls}}^* ,\lambda _{\mathrm{li}}^*$.
By Facts \ref{T1212}(a) and \ref{T1261}(a),   $\lambda _{\mathrm{s}}$ is a Hausdorff convergence, while
$\lambda _{\mathrm{ls}} $ and $ \lambda_{\mathrm{li}}$ are not; thus, $\lambda _{\mathrm{s}}< \lambda _{\mathrm{ls}},\lambda_{\mathrm{li}}$.
By Fact \ref{T1220}(b) $\lambda _{\mathrm{s}}^*$ is a Hausdorff convergence and, clearly, $\lambda _{\mathrm{ls}}^* $ and $ \lambda_{\mathrm{li}}^*$ are not
Hausdorff; so, $\lambda _{\mathrm{s}}^* < \lambda _{\mathrm{ls}}^* , \lambda _{\mathrm{li}}^*$.

By the construction of the closure $\lambda ^*$ it follows that we always have $\lambda \leq \lambda ^*$; thus
$\lambda _{\mathrm{ls}}\leq \lambda _{\mathrm{ls}}^*$,
$\lambda _{\mathrm{li}}\leq \lambda _{\mathrm{li}}^*$ and
$\lambda _{\mathrm{s}}\leq \lambda _{\mathrm{s}}^*$.
By Fact \ref{T1261}(f) we have
$ \lambda _{\mathrm{ls}}^* \leq \lim_{{\mathcal O}_{\lambda _{\mathrm{ls}}}}$ and
$\lambda _{\mathrm{li}}^* \leq \lim_{{\mathcal O}_{\lambda _{\mathrm{li}}}}$.
%$\lambda _{\mathrm{s}}^* \leq \lim_{{\mathcal O}_{\lambda _{\mathrm{s}}}}$
%follow from the fact that each topological convergence satisfies (L3).
%
The equality $\lambda _{\mathrm{s}}^* = \lim_{{\mathcal O}_{\lambda _{\mathrm{s}}}}$ follows from Facts \ref{T1212}(a) and \ref{T1220}(a).
Since ${\mathcal O}_{\mathrm{lsi}} \subset {\mathcal O}_{\lambda _{\mathrm{s}}}$ we have
$\lim_{{\mathcal O}_{\lambda _{\mathrm{s}}}} \leq \lim_{{\mathcal O}_{\mathrm{lsi}}}$.

Further we prove that
$\lim_{{\mathcal O}_{\mathrm{lsi}}}  = \lim_{{\mathcal O}_{\lambda_{\mathrm{ls} }}} \cap \lim_{{\mathcal O}_{\lambda_{\mathrm{li} }}}$.
Since ${\mathcal O}_{\lambda_{\mathrm{ls} }}, {\mathcal O}_{\lambda_{\mathrm{li} }} \subset {\mathcal O}_{\mathrm{lsi}}$,
%by Fact \ref{T1200}
we have $\textstyle \lim_{{\mathcal O}_{\mathrm{lsi}}}\leq \lim_{{\mathcal O}_{\lambda_{\mathrm{ls} }}}, \lim_{{\mathcal O}_{\lambda_{\mathrm{li} }}}$.
Conversely, if $a\in \lim_{{\mathcal O}_{\lambda_{\mathrm{ls} }}}(x) \cap \lim_{{\mathcal O}_{\lambda_{\mathrm{li} }}}(x)$
and $U$ is a ${\mathcal O}_{\mathrm{lsi}}$-neighborhood of $a$, then there is $O_1\cap O_2\in {\mathcal B}_{\mathrm{lsi}}$
such that $a\in O_1 \cap O_2\subset U$ and, hence, there are
$n_i \in \omega$, $i\in \{ 1,2 \}$, such that $x_n \in  O_i$, for
each $n\geq n_i$. Thus for each $n\geq \max \{ n_1, n_2 \}$ we have
$x_n \in U$, so $a\in \lim_{{\mathcal O}_{\mathrm{lsi}}}(x)$.

So we have $\lim_{{\mathcal O}_{\mathrm{lsi}}}  \leq \lim_{{\mathcal O}_{\lambda_{\mathrm{ls} }}} , \lim_{{\mathcal O}_{\lambda_{\mathrm{li} }}}$.
Assuming that $\lim_{{\mathcal O}_{\mathrm{lsi}}} = \lim_{{\mathcal O}_{\lambda_{\mathrm{ls} }}}$, we would have
$\lambda_{\mathrm{ls}}\leq\lim_{{\mathcal O}_{\lambda_{\mathrm{ls} }}} \leq \lim_{{\mathcal O}_{\lambda_{\mathrm{li} }}}$
and, since $1\in \lambda_{\mathrm{ls}}(\langle 0\rangle)$, %we would have
$1\in \lim_{{\mathcal O}_{\lambda_{\mathrm{li} }}}(\langle 0\rangle)$.
Now, since the sets from  ${\mathcal O}_{\lambda_{\mathrm{li} }}$ are upward-closed, for a non-empty set $O\in {\mathcal O}_{\lambda_{\mathrm{li} }}$
we would have $1\in O$ and, since $1\in \lim_{{\mathcal O}_{\lambda_{\mathrm{li} }}}(\langle 0\rangle)$, $0\in O$ as well,
which would give $O={\mathbb B}$.  So ${\mathcal O}_{\lambda_{\mathrm{li} }}$ would be the antidiscrete topology
which is false, because it is $T_0$. Thus $\lim_{{\mathcal O}_{\mathrm{lsi}}} < \lim_{{\mathcal O}_{\lambda_{\mathrm{ls} }}}$
and, similarly, $\lim_{{\mathcal O}_{\mathrm{lsi}}} < \lim_{{\mathcal O}_{\lambda_{\mathrm{li} }}}$, which implies that
${\mathcal O}_{\mathrm{lsi}} \supsetneq {\mathcal O}_{\lambda_{\mathrm{ls} }} , {\mathcal O}_{\lambda_{\mathrm{li} }}$.
\kdok
In the sequel we consider the topology ${\mathcal O}_{\mathrm{lsi}}$ and its convergence
and investigate the form of the diagrams in Figure \ref{FIG1} for particular (classes of) Boolean algebras.
In particular, it is natural to ask for which complete Boolean algebras we have
\begin{equation}\label{EQ004}
{\mathcal O}_{\mathrm{lsi}} = {\mathcal O}_{\lambda _{\mathrm{s}} }\; \mbox{ or, at least, }
\;\textstyle \lim_{{\mathcal O}_{\mathrm{lsi}}}=\lim_{{\mathcal O}_{\lambda _{\mathrm{s}} }}?
\end{equation}
First we give some sufficient conditions for these equalities.
\begin{te}\label{T1239}
Let ${\mathbb B}$ be a complete Boolean algebra. Then

(a) If the algebra ${\mathbb B}$ satisfies condition $(\hbar)$,
then $\textstyle \lim_{{\mathcal O}_{\mathrm{lsi}}}=\lim_{{\mathcal O}_{\lambda _{\mathrm{s}} }}$;

(b) If the algebra ${\mathbb B}$ is $(\omega , 2)$-distributive,
then $\lim_{{\mathcal O}_{\lambda _{\mathrm{ls}} }}=\lambda _{\mathrm{ls}} $,
$\lim_{{\mathcal O}_{\lambda _{\mathrm{li}} }}=\lambda _{\mathrm{li}} $ and
$\lim_{{\mathcal O}_{\mathrm{lsi}}}=\lim_{{\mathcal O}_{\lambda _{\mathrm{s}} }}=\lambda _{\mathrm{s}}$;
so the diagram for convergences collapses to %$\lambda _{\mathrm{ls}}$, $\lambda _{\mathrm{ls}}$ and $\lambda _{\mathrm{ls}}$
3 nodes;

(c) If $\,\Lim_{{\mathcal O}_{\mathrm{lsi}}}=\Lim_{{\mathcal O}_{\lambda _{\mathrm{s}} }}$,
then ${\mathcal O}_{\mathrm{lsi}}={\mathcal O}_{\lambda _{\mathrm{s}} }$ iff
    $\;\langle {\mathbb B}, {\mathcal O}_{\mathrm{lsi}} \rangle$ is a sequential space.
\end{te}
\dok (a) By Theorem \ref{T1286} we have ${\mathcal O}_{\mathrm{lsi}}\subset {\mathcal O}_{\lambda _{\mathrm{s}} }$
so, $\lim_{{\mathcal O}_{\lambda _{\mathrm{s}} }}\leq \lim_{{\mathcal O}_{\mathrm{lsi}}}$.

Conversely, assuming that  $x\in {\mathbb B}^\omega$ and $a \in \lim_{{\mathcal O}_{\mathrm{lsi}}} (x)$, by Theorem \ref{T1286} we have
\begin{equation}\label{EQ003}\textstyle
a \in \lim_{{\mathcal O}_{\lambda _{\mathrm{ls} }}}(x)\cap \lim_{{\mathcal O}_{\lambda _{\mathrm{li} }}}(x)
\end{equation}
and we prove that $a\in \lim_{{\mathcal O}_{\lambda _{\mathrm{s}} }}(x)$.
Thus, by Facts \ref{T1212}(a) and \ref{T1220}(a), we have to show
that for each $y \prec x$ there is $z \prec y$ such that $\limsup z=\liminf z=a$.

Let $y$ be a subsequence of $x$.
By Fact \ref{T1261}(b) the convergence $\lambda _{\mathrm{ls}} $ is weakly topological
so, by (\ref{EQ003}) and Fact \ref{T1220}(a),
there is $z' \prec y$ such that $\limsup z' \leq a$.
Since $z' \prec x$ and the convergence $\lambda _{\mathrm{li}} $ is weakly topological,
by (\ref{EQ003}) and Fact \ref{T1220}(a) again,
there is $z \prec z'$ such that $\liminf z \geq a$.
Now, we have $\limsup z \leq \limsup z' \leq a \leq \liminf z$,
which implies that $\liminf z = \limsup z =a$.

(b) If the algebra ${\mathbb B}$ is $(\omega , 2)$-distributive,
then by Facts \ref{T1212}(b) and \ref{T1261}(c) we have
$\lim_{{\mathcal O}_{\lambda _{\mathrm{s}} }}=\lambda _{\mathrm{s}} $,
$\lim_{{\mathcal O}_{\lambda _{\mathrm{ls}} }}=\lambda _{\mathrm{ls}} $ and
$\lim_{{\mathcal O}_{\lambda _{\mathrm{li}} }}=\lambda _{\mathrm{li}} $.
Thus, by Theorem \ref{T1286} we have
$\lim_{{\mathcal O}_{\mathrm{lsi}}}
=\lim_{{\mathcal O}_{\lambda _{\mathrm{ls}} }} \cap \lim_{{\mathcal O}_{\lambda _{\mathrm{li}} }}
=\lambda _{\mathrm{ls}} \cap \lambda _{\mathrm{li}}
= \lambda _{\mathrm{s}}
=\lim_{{\mathcal O}_{\lambda _{\mathrm{s}} }}$.

(c) The implication ``$\Rightarrow$" is true because the topology ${\mathcal O}_{\lambda _{\mathrm{s}} }$ is sequential.
If $\Lim_{{\mathcal O}_{\mathrm{lsi}}}=\Lim_{{\mathcal O}_{\lambda _{\mathrm{s}} }}$, and ${\mathcal O}_{\mathrm{lsi}} $ is a sequential topology,
then (since a topology ${\mathcal O}$ is sequential iff  ${\mathcal O}={\mathcal O}_{\lim _{\mathcal O}}$)
we have
${\mathcal O}_{\mathrm{lsi}}
={\mathcal O}_{\Lim_{{\mathcal O}_{\mathrm{lsi}}}}
={\mathcal O}_{\Lim_{{\mathcal O}_{\lambda _{\mathrm{s}} }}}
={\mathcal O}_{\lambda _{\mathrm{s}} }$.
\hfill $\Box$
\paragraph{The unit interval}
Although the unit interval $I=[0,1]$ is not a Boolean algebra,
it provides obvious examples of the convergences considered in this paper.
Let ${\mathcal O}_{\leftarrow}=\{ [0,a) :0<a\leq 1\}\cup \{ \emptyset ,I\}$ and
${\mathcal O}_{\rightarrow}=\{ (a,1] :0\leq a < 1\}\cup \{ \emptyset ,I\}$ be the left and the right topology on $I$
and let ${\mathcal O}_{\mathrm{st}}$ denote the standard topology on $I$.
It is easy to check that
\begin{equation}\label{EQ000}\textstyle
\lim_{{\mathcal O}_{\leftarrow}}=\lambda _{\mathrm{ls}} \;\;\mbox{ and } \;\;
\lim_{{\mathcal O}_{\rightarrow}}=\lambda _{\mathrm{li}}\;\;\mbox{ and } \;\;
\lim_{{\mathcal O}_{\mathrm{st}}}=\lambda _{\mathrm{s}}.
\end{equation}
By Theorem 2.6 of \cite{KuPaNSJOM} a topology ${\mathcal O}$ is sequential iff ${\mathcal O}={\mathcal O}_{\lim _{\mathcal O}}$.
So since the topology ${\mathcal O}_{\leftarrow}$ is first countable and, hence, sequential, by (\ref{EQ000}) we have
${\mathcal O}_{\lambda_{\mathrm{ls} }}={\mathcal O}_{\lim_{{\mathcal O}_{\leftarrow}}}={\mathcal O}_{\leftarrow}$; and similarly for the other two topologies.
So
\begin{equation}\label{EQ001}\textstyle
{\mathcal O}_{\lambda_{\mathrm{ls} }}={\mathcal O}_{\leftarrow} \;\;\mbox{ and } \;\;
{\mathcal O}_{\lambda_{\mathrm{li} }}={\mathcal O}_{\rightarrow} \;\;\mbox{ and } \;\;
{\mathcal O}_{\lambda_{\mathrm{s} }}={\mathcal O}_{\mathrm{st}}.
\end{equation}
Since ${\mathcal O}_{\leftarrow}\cup {\mathcal O}_{\rightarrow}$ is a subbase of ${\mathcal O}_{\mathrm{st}}$
we have ${\mathcal O}_{\mathrm{st}}={\mathcal O}({\mathcal O}_{\leftarrow}\cup {\mathcal O}_{\rightarrow})$
and by (\ref{EQ001}) we have
${\mathcal O}_{\lambda_{\mathrm{s} }}
={\mathcal O}_{\mathrm{st}}
={\mathcal O}({\mathcal O}_{\leftarrow}\cup {\mathcal O}_{\rightarrow})
={\mathcal O}({\mathcal O}_{\lambda_{\mathrm{ls} }} \cup {\mathcal O}_{\lambda_{\mathrm{li} }})
={\mathcal O}_{\lambda_{\mathrm{lsi} }}
$
and (\ref{EQ004}) is true.
\paragraph{Power set algebras}
Let $\kappa \geq \omega$ be a cardinal. We recall that the {\it Alexandrov cube of weight $\kappa$}
%in notation $\mathbb{A}_\kappa =\langle 2^\kappa , \tau_A \rangle$,
is the product of $\kappa$ many copies of the two point space $2=\{0,1\}$ with the topology $\{\emptyset, \{0\},\{0,1\}\}$.
%(This is a universal $T_0$-space of weight $\leq \kappa$.)
Identifying the sets  $P(\kappa)$ and $2^\kappa$ via characteristic functions
we obtain a homeomorphic copy $\mathbb{A}_\kappa =\langle P(\kappa) , \tau_{\mathbb{A}_\kappa} \rangle$
of that space.\footnote{
We recall that for a sequence $\langle X_n :n \in \omega\rangle$ in $P(\kappa )$ we have

$\textstyle \liminf _{n\in \omega}X_n \; =  \textstyle \bigcup _{k\in \omega}\bigcap _{n\geq k}X_n
                           =  \{ x: x\in X_n \mbox{ for all but finitely many }n \}$,

$\textstyle \limsup _{n\in \omega}X_n  =  \textstyle \bigcap _{k\in \omega}\bigcup _{n\geq k}X_n
                          =  \{ x: x\in X_n \mbox{ for infinitely many }n \}$.
}
Further, the {\it Cantor cube of weight $\kappa$}
%in notation $\mathbb{A}_\kappa =\langle 2^\kappa , \tau_A \rangle$,
is the product of $\kappa$ many copies of the two point discrete space $2=\{0,1\}$ and,
%(This is a universal  zero-dimensional space of weight $\leq \kappa$.)
identifying the sets  $P(\kappa)$ and $2^\kappa$ again, we obtain its homeomorphic copy
$\mathbb{C}_\kappa =\langle P(\kappa) , \tau_{\mathbb{C}_\kappa} \rangle$.
By \cite{KuPaAlx} we have
\begin{fac} \label{T1260}
For the power algebra $P(\kappa )$ with the Aleksandrov topology we have

(a) $\lambda _{\mathrm{ls} } = \lim _{{\mathcal O}_{\lambda _{\mathrm{ls} }}} = \lim_{\tau _{\mathbb{A}_\kappa}}$;
    thus $\lambda _{\mathrm{ls} } $ is a topological convergence;

(b) $\langle P(\kappa ), \tau_{\mathbb{A}_\kappa} \rangle $ is a sequential space
    iff ${\mathcal O}_{\lambda _{\mathrm{ls} }}= \tau _{\mathbb{A}_\kappa}$
     iff $\kappa = \omega$;

(c) If $\kappa > \omega$, then $\tau _{\mathbb{A}_\kappa}\varsubsetneq {\mathcal O}_{\lambda _{\mathrm{ls} }}\not\subset \tau _{\mathbb{C}_\kappa}$.

\noindent
For the power algebra $P(\kappa )$ with the Cantor topology we have

(d) $\lambda _{\mathrm{s} }= \lim _{{\mathcal O}_{\lambda _{\mathrm{s} }}} = \lim_{\tau _{\mathbb{C}_\kappa}}$;
    thus $\lambda _{\mathrm{s} } $ is a topological convergence;

(e) $\langle P(\kappa ), \tau_{\mathbb{C}_\kappa} \rangle $ is a sequential space
    iff ${\mathcal O}_{\lambda _{\mathrm{s} }}= \tau _{\mathbb{C}_\kappa}$
    iff $\kappa = \omega$;

(f) If $\kappa > \omega$, then $\tau _{\mathbb{C}_\kappa} \varsubsetneq {\mathcal O}_{\lambda _{\mathrm{s}} }$.
\end{fac}
Let $\tau _{\mathbb{A}_\kappa ^c}$  be the topology on the power algebra $P(\kappa )$ obtained by the standard identification
of $P(\kappa )$ and $2^\kappa$ with the
Tychonov  topology of $\kappa$ many copies
of the space 2 with the topology $\{\emptyset,\{1 \}, \{0,1\}\}$. Then, clearly, $X\mapsto \kappa \setminus X$ is a homeomorphism
from $\mathbb{A}_\kappa =\langle P(\kappa) , \tau_{\mathbb{A}_\kappa} \rangle$ onto the reversed Alexandrov cube
$\mathbb{A}_\kappa ^c =\langle P(\kappa) , \tau_{\mathbb{A}_\kappa^c} \rangle$.
Replacing $\tau_{\mathbb{A}_\kappa}$ by $\tau_{\mathbb{A}_\kappa^c}$
and $\lambda _{\mathrm{ls} }$ by $\lambda _{\mathrm{li} }$ in (a), (b) and (c) of Fact \ref{T1260} we obtain the corresponding dual statements.
In addition, we have
\begin{te}\label{T1235}
For the power algebra $P(\kappa )$ we have

(a) $\lim _{{\mathcal O}_{\mathrm{lsi}}}=\lambda _{\mathrm{s} }$;

(b) $\tau _{\mathbb{C}_\kappa}$ is the minimal topology containing $\tau_{\mathbb{A}_\kappa}$ and $\tau_{\mathbb{A}_\kappa^c}$;

(c) $\tau _{\mathbb{C}_\kappa}\subset {\mathcal O}_{\mathrm{lsi}}$ and
    ${\mathcal O}_{\mathrm{lsi}}$ is a Hausdorff topology on $P(\kappa)$;

(d) For $\kappa =\omega $ we have ${\mathcal O}_{\mathrm{lsi}}= \tau _{\mathbb{C}_\omega } = {\mathcal O}_{\lambda _{\mathrm{s}} }$;

(e) ${\mathcal O}_{\mathrm{lsi}}= \tau _{\mathbb{C}_\kappa}$ iff $\kappa =\omega $.
\end{te}
\dok
(a) By Fact \ref{T1260}(a) and its dual we have
$\lim _{{\mathcal O}_{\lambda _{\mathrm{ls} }}}=\lambda _{\mathrm{ls} }$ and $\lim _{{\mathcal O}_{\lambda _{\mathrm{li} }}}=\lambda _{\mathrm{li} }$.
Now, by Theorem \ref{T1286},
$\lim _{{\mathcal O}_{\mathrm{lsi}}}
=\lim _{{\mathcal O}_{\lambda _{\mathrm{ls} }}} \cap \lim _{{\mathcal O}_{\lambda _{\mathrm{li} }}}
=\lambda _{\mathrm{ls} } \cap \lambda _{\mathrm{li} }
=\lambda _{\mathrm{s} }$.

(b) Let ${\mathcal O}$ be the minimal topology containing $\tau_{\mathbb{A}_\kappa}$ and $\tau_{\mathbb{A}_\kappa^c}$.
A subbase for the topology $\tau_{\mathbb{A}_\kappa}$ (resp.\ $\tau_{\mathbb{A}_\kappa^c}$) consists of the sets
$B_i := \{ X\subset \kappa :i\not\in X\}$ (resp.\ $B_i^c := \{ X\subset \kappa :i\in X\}$), where $i\in \kappa$; while
the family $\mathcal{S}_{\mathbb{C}_\kappa}:=\bigcup _{i\in \kappa}\{B_i , B_i^c\}$ is a subbase for the topology $\tau _{\mathbb{C}_\kappa}$.
Thus $\tau_{\mathbb{A}_\kappa}\cup \tau_{\mathbb{A}_\kappa^c}\subset \tau _{\mathbb{C}_\kappa}$
and, hence, ${\mathcal O}\subset \tau _{\mathbb{C}_\kappa}$.
On the other hand, $\mathcal{S}_{\mathbb{C}_\kappa}\subset \tau_{\mathbb{A}_\kappa}\cup \tau_{\mathbb{A}_\kappa^c}\subset {\mathcal O}$, which gives
$\tau _{\mathbb{C}_\kappa}\subset {\mathcal O}$.

(c) By Fact \ref{T1260} and its dual we have $\tau_{\mathbb{A}_\kappa}\subset {\mathcal O}_{\lambda_{\mathrm{ls}}}$
and $\tau_{\mathbb{A}_\kappa^c} \subset {\mathcal O}_{\lambda_{\mathrm{li}}}$.
Thus  $\tau_{\mathbb{A}_\kappa}\cup \tau_{\mathbb{A}_\kappa^c} \subset {\mathcal O}_{\mathrm{lsi}} $ and
$\tau _{\mathbb{C}_\kappa} \subset {\mathcal O}_{\mathrm{lsi}}$, by the minimality of $\tau _{\mathbb{C}_\kappa}$ proved in (b).

(d) By (c) and Theorem \ref{T1286},
$\tau _{\mathbb{C}_\omega} \subset {\mathcal O}_{\mathrm{lsi}} \subset {\mathcal O}_{\lambda _{\mathrm{s}} }$
and we apply Fact \ref{T1260}(e).

(e) By (d) the implication ``$\Leftarrow$" is true.
Assuming that  ${\mathcal O}_{\mathrm{lsi}}= \tau _{\mathbb{C}_\kappa}$ and $\kappa >\omega $,
by Fact \ref{T1260}(c) we would have ${\mathcal O}_{\lambda _{\mathrm{ls} }}\not\subset \tau _{\mathbb{C}_\kappa}$, which gives a contradiction
because ${\mathcal O}_{\lambda _{\mathrm{ls} }}\subset{\mathcal O}_{\mathrm{lsi}}$.
\kdok

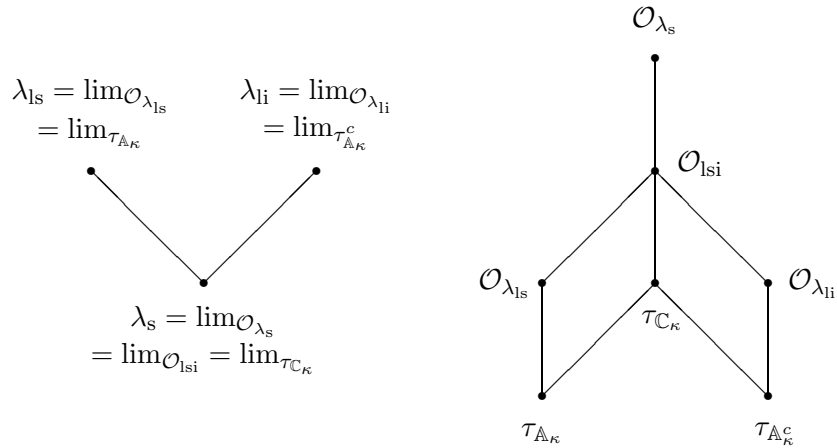
\begin{figure}[h]
\begin{center}
      \unitlength 1mm % = .854pt
\linethickness{0.4pt}
\ifx\plotpoint\undefined\newsavebox{\plotpoint}\fi % GNUPLOT compatibility
\begin{picture}(120,65)(0,0)
%====================== linije ===========================
\put(30,25){\line(-1,1){15}}%1
\put(30,25){\line(1,1){15}}%2
\put(75,10){\line(1,1){15}}%3
\put(105,10){\line(-1,1){15}}%4
\put(75,10){\line(0,1){15}}%5
\put(90,25){\line(0,1){15}}%6
\put(105,10){\line(0,1){15}}%7
\put(75,25){\line(1,1){15}}%8
\put(105,25){\line(-1,1){15}}%9
\put(90,40){\line(0,1){15}}%10
%====================== tacke ===========================
%\put(150,190){\circle*{2}}
\put(30,25){\circle*{1}}%1
\put(15,40){\circle*{1}}%2
\put(45,40){\circle*{1}}%3
\put(75,10){\circle*{1}}%4
\put(105,10){\circle*{1}}%5
\put(75,25){\circle*{1}}%6
\put(90,25){\circle*{1}}%7
\put(105,25){\circle*{1}}%8
\put(90,40){\circle*{1}}%9
\put(90,55){\circle*{1}}%10
%====================== tekst ===========================
%\put(165,190){\makebox(0,0)[cc]{${\mathcal O}_{\lambda_{\mathrm{li}}}$}}
\put(30,20){\makebox(0,0)[cc]{$ \lambda _{\mathrm{s} }= \lim _{{\mathcal O}_{\lambda _{\mathrm{s} }}}  $}}%1
\put(30,15){\makebox(0,0)[cc]{$ =\lim _{{\mathcal O}_{\mathrm{lsi}}} = \lim_{\tau _{\mathbb{C}_\kappa}} $}}%2
\put(15,50){\makebox(0,0)[cc]{$ \lambda _{\mathrm{ls} } = \lim _{{\mathcal O}_{\lambda _{\mathrm{ls} }}}  $}}%3
\put(15,45){\makebox(0,0)[cc]{$ = \lim_{\tau _{\mathbb{A}_\kappa}}  $}}%4
\put(45,50){\makebox(0,0)[cc]{$ \lambda _{\mathrm{li} } = \lim _{{\mathcal O}_{\lambda _{\mathrm{li} }}}  $}}%5
\put(45,45){\makebox(0,0)[cc]{$ = \lim_{\tau _{\mathbb{A}_\kappa}^c}  $}}%6
\put(75,5){\makebox(0,0)[cc]{$ \tau _{\mathbb{A}_\kappa}  $}}%7
\put(106,5){\makebox(0,0)[cc]{$  \tau _{\mathbb{A}_\kappa ^c}   $}}%8
\put(91,20){\makebox(0,0)[cc]{$ \tau _{\mathbb{C}_\kappa}  $}}%9
\put(70,25){\makebox(0,0)[cc]{$ {\mathcal O}_{\lambda _{\mathrm{ls} }}  $}}%10
\put(111,25){\makebox(0,0)[cc]{$ {\mathcal O}_{\lambda _{\mathrm{li} }}    $}}%11
\put(96,41){\makebox(0,0)[cc]{$ {\mathcal O}_{\mathrm{lsi}}  $}}%12
\put(90,60){\makebox(0,0)[cc]{$ {\mathcal O}_{\lambda _{\mathrm{s} }}   $}}%13
\end{picture}
\end{center}
\vspace{-8mm}
  \caption{Convergences and topologies on the algebra $P(\kappa )$}\label{FIG2}
\end{figure}
\noindent
For the power set algebras the diagrams from Figure \ref{FIG1} are presented in Figure \ref{FIG2}.
Namely, by Theorem \ref{T1239}(b), the diagram describing convergences collapses to
three nodes.
The diagram for topologies in Figure \ref{FIG2} contains the topologies from Figure \ref{FIG1} as well as the topologies of the
Cantor, Alexandrov and reversed Alexandrov cube (see Fact \ref{T1260}(c) and Theorem \ref{T1235}(c)).
By Fact \ref{T1260}(b) and (e), for $\kappa = \omega$
the diagram  describing topologies contains exactly three different topologies. So, for the algebra $P( \omega)$ we have
${\mathcal O}_{\mathrm{lsi}} ={\mathcal O}_{\lambda _{\mathrm{s} }}$ and (\ref{EQ004}) is true.
\paragraph{Maharam algebras}
We recall that a {\it submeasure} on a complete
Boolean algebra ${\mathbb B}$ is a function $\mu:{\mathbb B}\rightarrow[0,\infty)$ satisfying: (i) $\mu(0)=0$;
(ii) $a\leq b\Rightarrow \mu(a)\leq\mu(b)$ and (iii) $\mu(a\lor b)\leq\mu(a)+\mu(b)$.
A submeasure $\mu$ is {\it strictly positive} iff (iv) $a>0\Rightarrow \mu(a)>0$; $\mu$ is called a {\it Maharam}
(or a {\it continuous}) {\it submeasure} iff (v) $\lim_{n\rightarrow\infty}\mu(a_n)=0$ holds
for each decreasing sequence $\langle a_n : n\in\omega \rangle$ in ${\mathbb B}$
satisfying $\bigwedge_{n\in\omega }a_n=0$.
Then $\lim_{n\rightarrow\infty}\mu(a_n)=\mu (\bigwedge _{n\in \omega} a_n)$,
for each decreasing sequence $\langle a_n \rangle$ in ${\mathbb B}$.
A complete Boolean algebra ${\mathbb B}$ admitting a strictly positive Maharam submeasure is called a {\it Maharam algebra}.
\begin{te}  \label{T1262}
On each  Maharam algebra ${\mathbb B}$ we have ${\mathcal O}_{\mathrm{lsi}} = {\mathcal O}_{\lambda _{\mathrm{s}} }$.
\end{te}
\dok
Under the assumption, $d(x,y)=\mu (x\vartriangle y)$ is a metric on ${\mathbb B}$ which generates the topology
${\mathcal O}_{\lambda _{\mathrm{s}} }$ (see \cite{Mah47}). For a non-empty set $O\in {\mathcal O}_{\lambda _{\mathrm{s}} }$ we show that
$O\in {\mathcal O}_{\mathrm{lsi}}$. Let $a\in O$ and $r>0$, where
$B(a,r)= \{ x\in {\mathbb B} : \mu (x\vartriangle a)<r \} \subset O$. Let
$$
O_1 = \{ x\in {\mathbb B} : \mu (x\setminus a)<r/2 \} \mbox{ and }
O_2 = \{ x\in {\mathbb B} : \mu (a\setminus x)<r/2 \}.
$$
Then by (i) we have $a\in O_1 \cap O_2$. If $x\in O_1 \cap O_2$, then, by (iii),
$\mu (x\vartriangle a) \leq \mu (x\setminus a)+ \mu (a\setminus x)<r$ and, hence,
$x\in B(a,r)$, thus $a\in O_1 \cap O_2 \subset O$.

Let us prove that $O_1\in {\mathcal O}_{\lambda _{{\mathrm{ls} }}}$.
By Fact \ref{T1261}(a) the convergence $\lambda _{{\mathrm{ls}}}$ satisfies (L1) and (L2), so  it is
sufficient to prove that ${\mathbb B}\setminus O_1$ is a closed set,
which means that $u_{\lambda _{{\mathrm{ls} }}}({\mathbb B}\setminus O_1)\subset {\mathbb B}\setminus O_1$.
By (iii), the set
${\mathbb B}\setminus O_1$ is upward-closed and it is sufficient to
show that for a sequence $\langle x_n \rangle$ in ${\mathbb
B}\setminus O_1$ we have $\limsup x_n \in {\mathbb B}\setminus O_1$,
that is $\mu (\limsup x_n \setminus a )\geq r/2$. By the assumption
we have $\mu (x_n \setminus a )\geq r/2$, for each $n\in \omega$.
Now $\limsup x_n \setminus a = \bigwedge _{k\in \omega  }y_k$, where
$y_k = \bigvee _{n\geq k} x_n \setminus a$, $k\in \omega $, is a
decreasing sequence and $\mu (y_k)\geq r/2$ so, by the continuity of
$\mu$, $\mu (\limsup x_n \setminus a )= \lim _{k\rightarrow \infty }
\mu (y_k)\geq r/2$. Similarly we prove that $O_2\in {\mathcal
O}_{\lambda _{{\mathrm{li} }}}$ so  $O_1 \cap O_2 \in {\mathcal
O}_{\mathrm{lsi}}$ and $O$ is an ${\mathcal O}_{\mathrm{lsi}}$-neighborhood of the point $a$.
\kdok

\begin{figure}[h]
\begin{center}
      \unitlength 1mm % = .854pt
\linethickness{0.4pt}
\ifx\plotpoint\undefined\newsavebox{\plotpoint}\fi % GNUPLOT compatibility
\begin{picture}(120,50)(0,0)
%====================== linije ===========================
\put(30,10){\line(-1,1){15}}%1
\put(30,10){\line(0,1){15}}%2
\put(30,10){\line(1,1){15}}%3
\put(15,25){\line(0,1){15}}%4
\put(30,25){\line(-1,1){15}}%5
\put(30,25){\line(1,1){15}}%6
\put(45,25){\line(0,1){15}}%7
\put(75,25){\line(1,1){15}}%8
\put(105,25){\line(-1,1){15}}%9
%====================== tacke ===========================
\put(30,10){\circle*{1}}%1
\put(15,25){\circle*{1}}%2
\put(30,25){\circle*{1}}%3
\put(45,25){\circle*{1}}%4
\put(15,40){\circle*{1}}%5
\put(45,40){\circle*{1}}%6
\put(75,25){\circle*{1}}%7
\put(105,25){\circle*{1}}%8
\put(90,40){\circle*{1}}%9
%====================== tekst ===========================
\put(30,5){\makebox(0,0)[cc]{$ \lambda _{\mathrm{s} }  $}}%1
\put(10,25){\makebox(0,0)[cc]{$ \lambda _{\mathrm{ls} }  $}}%2
%\put(31,25){\makebox(0,0)[cc]{$ =\;\;\lambda _{\mathrm{ls} }^*   $}}%3===========
%\put(30,33){\makebox(0,0)[cc]{$ \lim _{{\mathcal O}_{\mathrm{lsi} } }\! =$}}%3,5
\put(29,25){\makebox(0,0)[cc]{$ \lim _{{\mathcal O}_{\lambda _{\mathrm{ls} }}} \;= \lambda _{\mathrm{ls} }^* $}}%3
\put(50,25){\makebox(0,0)[cc]{$ \lambda _{\mathrm{li} }  $}}%4
\put(15,45){\makebox(0,0)[cc]{$ \lambda _{\mathrm{ls} }^* = \lim _{{\mathcal O}_{\lambda _{\mathrm{ls} }}}  $}}%5
\put(48,45){\makebox(0,0)[cc]{$ \lambda _{\mathrm{li} }^* = \lim _{{\mathcal O}_{\lambda _{\mathrm{li} }}}   $}}%6
\put(75,20){\makebox(0,0)[cc]{$  {\mathcal O}_{\lambda _{\mathrm{ls} }} $}}%7
\put(105,20){\makebox(0,0)[cc]{$ {\mathcal O}_{\lambda _{\mathrm{li} }}  $}}%8
\put(90,45){\makebox(0,0)[cc]{$  {\mathcal O}_{\lambda _{\mathrm{s} }}= {\mathcal O}_{\mathrm{lsi}} $}}%9
\end{picture}
\end{center}
\vspace{-10mm}
  \caption{Non-$(\omega ,2)$-distributive Maharam algebras}\label{FIG3}
\end{figure}
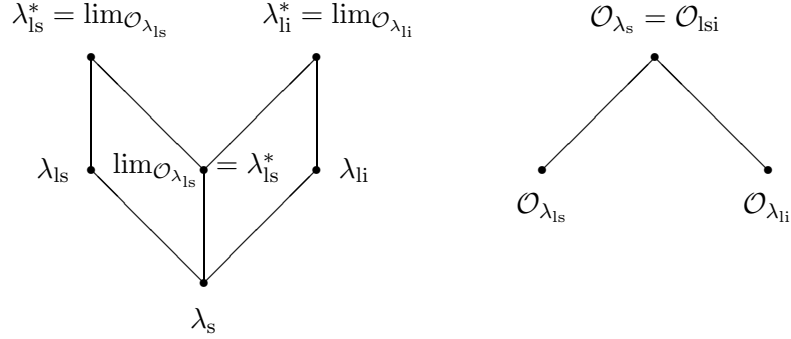
\noindent
Thus, if ${\mathbb B}$ is a Maharam algebra which is not $(\omega ,2)$-distributive
(for example, the algebra of the Lebesgue-measurable subsets of $[0,1]$ modulo the ideal of the sets of measure zero),
then, the Figure \ref{FIG3} describes the corresponding diagrams.
Namely, by Facts \ref{T1212}(a) and \ref{T1220}(a) we have $\lim_{\mathcal{O}_{\lambda _{\mathrm{s}}}}=\lambda _{\mathrm{s}}^*$
and, by Fact \ref{T1212}(b), $\lambda _{\mathrm{s}}\neq \lim_{\mathcal{O}_{\lambda _{\mathrm{s}}}}$.
Since the algebras with strictly positive measure satisfy the countable chain condition %ccc
the algebra ${\mathbb B}$ has $(\hbar )$.
Thus, by Facts \ref{T1261}(b) and \ref{T1220}(a) we have $\lim_{\mathcal{O}_{\lambda _{\mathrm{ls}}}}=\lambda _{\mathrm{ls}}^*$ and
$\lim_{\mathcal{O}_{\lambda _{\mathrm{li}}}}=\lambda _{\mathrm{li}}^*$.
By Fact \ref{T1261}(c) we have $\lambda _{\mathrm{ls}}\neq \lim_{\mathcal{O}_{\lambda _{\mathrm{ls}}}}$ and
$\lambda _{\mathrm{li}}\neq \lim_{\mathcal{O}_{\lambda _{\mathrm{li}}}}$.
By Theorem \ref{T1262} we have ${\mathcal O}_{\mathrm{lsi}} = {\mathcal O}_{\lambda _{\mathrm{s}} }$ and, hence,
$\lim_{{\mathcal O}_{\mathrm{lsi}}} =\lim _{ {\mathcal O}_{\lambda _{\mathrm{s}} }}$.
\paragraph{Collapsing algebras}
We show that both equalities from (\ref{EQ004}) can fail.
We recall that a family $T\subset [\omega]^\omega$ is a
{\it tower} iff it is well-ordered by $^*\!\!\supsetneq$ and has no pseudointersection.
The {\it tower number}, $\mathfrak{t}$, is the minimal cardinality of a tower.
A family ${\mathcal T}\subset [\omega ]^{\omega }$ is called a {\it base matrix tree} iff
$\langle {\mathcal T}, {}^* \!\! \supset \rangle$ is a tree of height ${\mathfrak h}$
and ${\mathcal T}$ is a dense set in the pre-order $\langle [\omega ]^{\omega }, \subset ^* \rangle $.
By a theorem of Balcar, Pelant and Simon (see \cite{Balc0}), such a tree always exists,
its levels are maximal almost disjoint families
and maximal chains in ${\mathcal T}$ are towers.
\begin{te}\label{T1289}
If  ${\mathbb B}$ is a complete Boolean algebra satisfying
$1\Vdash_{\mathbb B} ({\mathfrak h}^V)^{\check{~}} < {\mathfrak t}$ and {\rm cc}$({\mathbb B})>2^{\mathfrak h}$,
then  $\lim_{{\mathcal O}_{\lambda _{\mathrm{s}} }}<\lim_{{\mathcal O}_{\mathrm{lsi}}}$ and
${\mathcal O}_{\mathrm{lsi}}\varsubsetneq {\mathcal O}_{\lambda _{\mathrm{s}} }$.
\end{te}
\dok
Using the construction from the proof of Theorem 6.4 from \cite{KuPaAlx},
we will find a  sequence $x$ in ${\mathbb B}$ such that $0\in \lim_{{\mathcal O}_{\mathrm{lsi}}}(x) \setminus \lim_{{\mathcal O}_{\lambda _{\mathrm{s}} }}(x)$.

Let  ${\mathcal T}\subset [\omega]^\omega$ be a base matrix tree and $\Br({\mathcal T})$ the set of its maximal branches.
Since the height of ${\mathcal T}$ is ${\mathfrak h}$, the branches of ${\mathcal T}$ are of size $\leq {\mathfrak h}$;
so $\kappa :=| \Br({\mathcal T})| \leq {\mathfrak c}^{\mathfrak h}=2 ^{\mathfrak h}$
and we take a one-to-one  enumeration $\Br({\mathcal T})= \{ T_\alpha : \alpha < \kappa \}$.

Since $1\Vdash ({\mathfrak h}^V)^{\check{~}} <{\mathfrak t}$, for each $\alpha < \kappa $ we have
$1\Vdash |\check{T_\alpha }|<{\mathfrak t}$ and, consequently, in each generic extension of the ground model by ${\mathbb B}$
the family $T_\alpha$ obtains a pseudointersection.
Thus $1\Vdash \exists X \in [\check{\omega}]^{\check{\omega}} \; \forall B \in \check{T_\alpha }\; X \subset^*B$
so, by the Maximum Principle (see \cite[p.\ 226]{Kun}), there is
a name $\sigma_\alpha \in V^{\mathbb B}$ such that
\begin{equation}\label{EQ1290a}
1 \Vdash \sigma_\alpha  \in [\check{\omega}]^{\check{\omega}} \wedge
\forall B \in \check{T_{\alpha }} \;\; \sigma_\alpha \subset^*B.
\end{equation}
Since {\rm cc}$({\mathbb B})>2^{\mathfrak h}\geq \kappa$, there is a maximal antichain in  ${\mathbb B}$ of cardinality $\kappa$,
say $\{b_\alpha: \alpha < \kappa \}$.
By the Mixing lemma (see \cite[p.\ 226]{Kun})
there is a name $\tau \in V^{\mathbb B}$ such that
\begin{equation}\label{EQ1290b}
\forall \alpha <\kappa \;\; b_\alpha\Vdash \tau=\sigma_\alpha,
\end{equation}
and, clearly,  $1\Vdash \tau\in [\check{\omega}]^{\check{\omega}}$.
Let $x =\langle x_n \rangle \in {\mathbb B}^\omega$, where $x_n:=\|\check{n}\in \tau\|$, for $n \in \omega$.
Then for the corresponding name  $\tau _x = \{ \langle \check{n}, x_n \rangle : n\in \omega \}$ we have
\begin{equation}\label{EQ1290c}
1\Vdash \tau=\tau_x.
\end{equation}
Now, by Claims 1 and 2 from the proof of Theorem 6.4 from \cite{KuPaAlx} we have
$$
0  \in \Lim_{{\mathcal O}_{\lambda_{{\mathrm{ls}}}}}(x)\setminus \lambda^*_{{\mathrm{ls}}}(x).
$$
By Facts \ref{T1212}(a) and \ref{T1261}(g) we have
$\lim_{{\mathcal O}_{\lambda _{\mathrm{s}} }}(x)
= \lambda _{\mathrm{s}}^*(x)
=\lambda _{\mathrm{ls}}^*(x) \cap \lambda _{\mathrm{li}}^*(x)$
and, since $0\not\in \lambda _{\mathrm{ls}}^*(x)$, it follows that
$0\not\in \lim_{{\mathcal O}_{\lambda _{\mathrm{s}} }}(x)$.

By Theorem \ref{T1286} we have
$ \lim_{{\mathcal O}_{\mathrm{lsi}}}(x)
= \lim_{{\mathcal O}_{\lambda _{\mathrm{ls}} }}(x)\cap \lim_{{\mathcal O}_{\lambda _{\mathrm{li}} }}(x)$
and, since $0\in \lim_{{\mathcal O}_{\lambda _{\mathrm{ls}} }}(x)$, it remains to be proved that
$0\in \lim_{{\mathcal O}_{\lambda _{\mathrm{li}} }}(x)$.
But, if $0\in O \in {\mathcal O}_{\lambda _{\mathrm{li}}}$,
then, since $O$ is an upward-closed set, we have $O={\mathbb B}$. Consequently,
$x_n \in O$, for all $n\in \omega$, so $0\in \lim_{{\mathcal O}_{\lambda _{\mathrm{li}} }}(x)$.
\hfill $\Box$
\begin{ex}\rm\label{EX1209}
An algebra for which the diagrams describing convergences and topologies from Figure \ref{FIG1}
contain exactly 9 and 4 different objects respectively.

If ${\mathbb B}$ is a complete Boolean algebra which collapses $2^{\mathfrak h}$ to $\omega$
(e.g.\ the collapsing algebra $\Coll (\omega , 2^{\mathfrak h})=\mathrm{r.o.}({}^{<\omega }(2^{\mathfrak h}))$),
then ${\mathbb B}$ satisfies the assumptions of Theorem \ref{T1289} and, hence,
$\lim_{{\mathcal O}_{\lambda _{\mathrm{s}} }}<\lim_{{\mathcal O}_{\mathrm{lsi}}}$ and
${\mathcal O}_{\mathrm{lsi}}\varsubsetneq {\mathcal O}_{\lambda _{\mathrm{s}} }$.
By Theorem 6.4 from \cite{KuPaAlx} the same conditions provide that
the convergence  $\lambda _{\mathrm{ls}}$ is not weakly topological,
which, by Fact \ref{T1220}(a), gives $\lambda _{\mathrm{ls}}^*<  \lim _{{\mathcal O}_{\lambda _{\mathrm{ls}}}}$.
By Theorem 4.4 from \cite{KuPaAlx}, the mapping
$h:\langle {\mathbb B},{\mathcal O_{\lambda_{\mathrm{ls}}}}\rangle \rightarrow \langle {\mathbb B},{\mathcal O_{\lambda_{\mathrm{li}}}}\rangle$
given by $h(b)=b'$, for each $b\in  {\mathbb B}$, is a homeomorphism,
so $\lambda _{\mathrm{li}}^*<  \lim _{{\mathcal O}_{\lambda _{\mathrm{li}}}}$ as well.
Assuming that $\lambda _{\mathrm{ls}}=\lambda _{\mathrm{ls}}^*$,
by duality we would have $\lambda _{\mathrm{li}}=\lambda _{\mathrm{li}}^*$
and, by Theorem \ref{T1286},
$\lim _{{\mathcal O}_{\lambda _{\mathrm{s}}}}
=\lambda _{\mathrm{s}}^*
=\lambda _{\mathrm{ls}}^* \cap \lambda _{\mathrm{li}}^*
=\lambda _{\mathrm{ls}} \cap \lambda _{\mathrm{li}}
=\lambda _{\mathrm{s}}$.
But this is not true since the algebra ${\mathbb B}$ is not $(\omega ,2)$-distributive.
Thus $\lambda _{\mathrm{ls}}<\lambda _{\mathrm{ls}}^*$ and, similarly, $\lambda _{\mathrm{li}}<\lambda _{\mathrm{li}}^*$.
By Fact \ref{T1212}(b) we have $\lambda _{\mathrm{s}}<\lim_{{\mathcal O}_{\lambda _{\mathrm{s}} }}$.
The rest follows from Theorem \ref{T1286}.
\end{ex}
\noindent
{\bf Acknowledgments}
This research was supported by the Ministry of Education and Science of the Republic of Serbia (Project 174006).

\end{document}